\documentclass[12pt,a4paper]{article}

\usepackage[left=1cm,right=1cm,top=2cm,bottom=2cm]{geometry}
\usepackage[dvips]{graphics}
\usepackage[dvips]{epsfig}
\usepackage{color}
\usepackage{hyperref}
\usepackage{tikz,authblk}

\newtheorem{lemma}{Lemma}
\newtheorem{theorem}[lemma]{Theorem}
\newtheorem{corollary}[lemma]{Corollary}
\newtheorem{proposition}[lemma]{Proposition}
\newtheorem{definition}{Definition}
\newtheorem{remark}{Remark}

\newcommand{\dimo}{\noindent \emph{Proof. }}
\newcommand{\qed}{\\ \rightline{$\Box$ \ \ \ \ \ \ \ \ \ \ \ \ \ \ \ }\\}

\usepackage{latexsym}
\usepackage{amsfonts}
\usepackage{amssymb}
\usepackage{amsmath}
\usepackage[english]{babel}

\begin{document}
\title{Combinatorial properties of the G-degree}

 \renewcommand{\Authfont}{\scshape\small}
 \renewcommand{\Affilfont}{\itshape\small}
 \renewcommand{\Authand}{ and }
\author[1] {Maria Rita Casali}
\author[2] {Luigi Grasselli}

\affil[1] {Department of Physics, Mathematics and Computer Science, University of Modena and Reggio Emilia, \ \ \ \ \ \ \ \ \ Via Campi 213 B - 41125 Modena (Italy), casali@unimore.it}

\affil[2] {Department of Sciences and Methods for Engineering, University of Modena and Reggio Emilia, \ \ \ \ \ \ \ \ \ \ \ \ \ \ \ \ \ \ \ \ \ \ \ \ \ \ \  Via Amendola 2, Pad. Morselli - 42122 Reggio Emilia (Italy), luigi.grasselli@unimore.it}

\maketitle

\abstract {A strong interaction is known to exist between edge-colored graphs (which encode PL pseudo-manifolds of arbitrary dimension) and random tensor models (as a possible approach to the study of Quantum Gravity). The key tool is the {\it G-degree} of the involved graphs, which drives the {\it $1/N$ expansion} in the tensor models context.
In the present paper - by making use of combinatorial properties concerning Hamiltonian decompositions of the complete graph - we prove that, in any even dimension $d\ge 4$, the G-degree of all bipartite graphs, as well as of all (bipartite or non-bipartite) graphs representing singular manifolds, is an integer multiple of $(d-1)!$. 
As a consequence, in even dimension,  the terms of  the  $1/N$ expansion corresponding to odd powers of $1/N$ are null in the complex context, and do not involve colored graphs representing singular manifolds in the real context. 
\noindent
In particular, in the 4-dimensional case, where the G-degree is shown to depend only on the regular genera with respect to an arbitrary pair of  ``associated" cyclic permutations, 
several results are obtained, relating the G-degree or the regular genus of 5-colored graphs and the Euler characteristic of the associated PL 4-manifolds.}
 \endabstract



\bigskip
\section{Introduction}
\label{intro}
It is well-known that regular edge-colored graphs may encode PL-pseudomanifolds, giving rise to a combinatorial representation theory ({\it crystallization theory}) for singular PL-manifolds of arbitrary dimension (see Section \ref{sec: Preliminaries}).

In the last decade, the strong interaction between the topology of edge-colored graphs and random tensor models has been deeply investigated, bringing insights in both research fields.

The colored tensor models theory arises as a possible approach to the study of Quantum Gravity: in some sense, its aim is to generalize to higher dimension the matrix models theory which, in dimension two, has shown to be quite useful at providing a framework for Quantum Gravity. The key generalization is the recovery of the so called {\it $1/N$ expansion} in the tensor models context. In matrix models, the $1/N$ expansion is driven by the genera of the surfaces represented by Feynman graphs; in the higher dimensional setting of tensor models the $1/N$ expansion is driven by the {\it G-degree} of these graphs (see Definition \ref{G-degree}), that equals the genus of the represented surface in dimension two.

If $(\mathbb{C}^N)^{\otimes d }$ denotes the $d$-tensor product of the $N$-dimensional complex space $\mathbb{C}^N,$ a {\it $(d+1)$-dimensional colored tensor model} is a formal partition function
\begin{equation*}
\mathcal{Z}[N,\{\alpha_{B}\}]:=\int_{\mbox{f}}\frac{dTd\overline{T}}{(2\pi)^{N^d}}\exp(-N^{d-1}\overline{T}\cdot T + \sum_{B}\alpha_BB(T,\overline{T})),
\end{equation*}
where $T$ belongs to $(\mathbb{C}^N)^{\otimes d },$ $\overline{T}$ to its dual and $B(T,\overline{T})$ are {\it trace invariants} obtained by contracting the indices of the components of $T$ and $\overline{T}.$
In this framework, colored graphs naturally arise as Feynman graphs encoding tensor trace invariants.

As shown in \cite{BGR}, the {\it free energy}  $\frac{1}{N^d}\log \mathcal{Z}[N,\{t_{B}\}]$ is the formal series
\begin{equation} \label{1/N expansion}
\frac{1}{N^d}\log \mathcal{Z}[N,\{t_{B}\}] = \sum_{\omega_G\ge 0}N^{-\frac{2}{(d-1)!}\omega_G}F_{\omega_G}[\{t_B\}]\in \mathbb{C}[[N^{-1}, \{t_{B}\}]],
\end{equation}
\noindent where the coefficients $F_{\omega_G}[\{t_B\}]$ are generating functions of connected bipartite $(d+1)$-colored graphs with fixed G-degree $\omega_G$.

The {\it  $1/N$ expansion} of formula \eqref{1/N expansion} describes the r\^ole of colored graphs (and of their G-degree $\omega_G$) within colored tensor models theory and explains the importance of trying to understand which are the manifolds and pseudomanifolds represented by $(d+1)$-colored graphs with a given G-degree.

A more detailed description of these relationships between Quantum Gravity via tensor models and topology of colored graphs may be found in \cite{BGR}, \cite{Gurau-Ryan}, \cite{Gurau-book}, \cite{Casali-Cristofori-Dartois-Grasselli}.

A parallel tensor models theory, involving {\it real} tensor variables $T\in (\mathbb{R}^N)^{\otimes d}$, has been developed, taking into account also non-bipartite colored graphs (see \cite{Witten}): this is why both bipartite and non-bipartite colored graphs will be considered within the paper.

\smallskip 

Section \ref{sec: Preliminaries} contains a quick review of crystallization theory, including the idea  of {\it regular embedding} of edge-colored graphs into surfaces, which is crucial for the definitions of {\it G-degree} and {\it regular genus} of graphs (Definition \ref{G-degree}).

In Section \ref{sec: Factorizations}, combinatorial properties concerning Hamiltonian decompositions of the complete graph allow to prove the main results of the paper.

\smallskip 

\begin{theorem} \label{Th.1}
For each  bipartite $(d+1)$-colored graph $(\Gamma,\gamma)$,  with $d$ even, $d \ge 4$, \ $$\omega_G (\Gamma) \equiv 0 \mod (d-1)!$$
\end{theorem}

\begin{theorem} \label{Th.2}
For each  (bipartite or non-bipartite) $(d+1)$-colored graph $(\Gamma,\gamma)$ representing a singular $d$-manifold, with $d$ even, $d\ge 4,$  \ $$\omega_G (\Gamma) \equiv 0 \mod (d-1)!$$
\end{theorem}

Note that the above results turn out to have specific importance in the tensor models framework. In fact Theorem \ref{Th.1} implies that, in the $d$-dimensional complex context, with $d$ even and $d\ge 4$,  the only non-null terms in the  $1/N$ expansion of formula \eqref{1/N expansion} are the ones corresponding to even (integer) powers of $1/N.$

On the other hand, Theorem \ref{Th.2} ensures that in the real tensor models framework, where also non-bipartite graphs are involved, the $1/N$ expansion contains colored graphs representing (orientable or non-orientable)  singular manifolds - and, in particular, closed manifolds - only in the terms corresponding to even (integer) powers of $1/N.$
Both Theorems extend to arbitrary even dimension a result proved in \cite[Corollary 23]{Casali-Cristofori-Dartois-Grasselli} for graphs representing singular $4$-manifolds.

\smallskip 

Section \ref{sec: n=4} is devoted to the 4-dimensional case: in this particular situation, the general results of  Section \ref{sec: Factorizations} allow to obtain interesting properties relating
the G-degree with the topology of the associated PL 4-manifolds.  In fact, the G-degree of a $5$-colored graph is shown to depend only on the regular genera with respect to an arbitrary pair of  ``associated" cyclic permutations (Proposition \ref{iper-reduced G-degree(d=4)}). This fact yields relations between these two genera and the Euler characteristic of  the associated PL 4-manifold (Proposition \ref{inequalities Euler characteristic} and Proposition \ref{sum/difference}); moreover, two interesting classes of crystallizations arise in a natural way,  whose intersection consists in the known class of semi-simple crystallizations, introduced in \cite{[Basak-Casali 2016]} (see Remark \ref{rem:intersection-classes}).

\section{Edge-colored graphs and G-degree}
\label{sec: Preliminaries}

A {\it singular $d$-manifold} is a compact connected $d$-dimensional polyhedron admitting a simplicial triangulation where the links of vertices are closed connected $(d-1)$-manifolds, while the link of any $h$-simplex, for each $h > 0$, is a PL $(d-h-1)$-sphere.  A vertex whose link is not a PL $(d-1)$-sphere is called {\it singular}.

\begin{remark}\label{correspondence-sing-boundary} 
{\rm The class of singular $d$-manifolds includes the class of closed $d$-manifolds: in fact, a closed $d$-manifold is a singular $d$-manifold without singular vertices. 
Moreover, if $N$ is a singular $d$-manifold, then a compact PL $d$-manifold $\check N$ is obtained by deleting small open neighbourhoods of its singular vertices.
Obviously, $N=\check N$ if and only if $N$ is a closed manifold; otherwise, $\check N$ has a non-empty boundary without spherical components.
Conversely, given a compact PL $d$-manifold $M$, a singular $d$-manifold $\widehat M$ can be obtained by capping off each component of $\partial M$ by a cone over it.

Note that, in virtue of the above correspondence,
a bijection is defined between singular $d$-manifolds and compact PL $d$-manifolds with no spherical boundary components.
}
\end{remark}

\begin{definition} \label{$d+1$-colored graph}
A $(d+1)$-colored graph ($d \ge 2$) is a pair $(\Gamma,\gamma)$, where $\Gamma=(V(\Gamma), E(\Gamma))$ is a regular $d+1$ valent multigraph  (i.e. multiple edges are allowed, while loops are forbidden)  and $\gamma: E(\Gamma) \rightarrow \Delta_d=\{0,\ldots, d\}$ is a map injective on adjacent edges, called {\it coloration}.
\end{definition}

\smallskip 

For every $\mathcal B\subseteq\Delta_d$ let $\Gamma_{\mathcal B}$ be the subgraph obtained from $(\Gamma, \gamma)$ by deleting all the edges colored by $\Delta_d - \mathcal B$. The connected components of $\Gamma_{\mathcal B}$ are called {\it ${\mathcal B}$-residues} or, if $\#\mathcal B = h$, {\it $h$-residues} of $\Gamma$; the symbol $g_{\mathcal B}$ denotes their number. In the following, if $\mathcal B =\{c_1,\ldots,c_h\}$, its complementary set $\Delta_d - \mathcal B$ will be denoted by $\hat c_1\ldots\hat c_h.$

\smallskip

\noindent Given a $(d+1)$-colored graph $(\Gamma, \gamma)$, a $d$-dimensional pseudocomplex $K(\Gamma)$ can be associated by the following rules:
\begin{itemize}
\item for each vertex of $\Gamma$, let us consider a $d$-simplex and label its vertices by the elements of $\Delta_d$;
\item for each pair of $c$-adjacent vertices of $\Gamma$ ($c\in\Delta_d$), let us glue the corresponding $d$-simplices along their $(d-1)$-dimensional faces opposite to the $c$-labeled vertices, so that equally labeled vertices are identified.
\end{itemize}

\smallskip

$|K(\Gamma)|$ turns out to be a {\it $d$-pseudomanifold}\footnote{In fact, $|K(\Gamma)|$ is a {\it quasi-manifold}: see  \cite{Gagliardi 1979}.}, which  is orientable if and only if $\Gamma$ is bipartite,  and $(\Gamma, \gamma)$ is said to {\it represent} it.

\smallskip 

Note that, by construction, $K(\Gamma)$ is endowed with a vertex-labeling by $\Delta_d$ that is injective on any simplex.
Moreover, a bijective correspondence exists between the $h$-residues of $\Gamma$ colored by any $\mathcal B\subseteq\Delta_d$ and the $(d-h)$-simplices of $K(\Gamma)$ whose vertices are labeled by $\Delta_d - \mathcal B$.

In particular, for any color $c\in\Delta_d$, each connected component of $\Gamma_{\hat c}$ is a $d$-colored graph representing  a pseudocomplex that is PL-homeomorphic to the link of a $c$-labeled vertex of $K(\Gamma)$ in its first barycentric subdivision.
As a consequence, $|K(\Gamma)|$ is a singular $d$-manifold (resp. a closed $d$-manifold) iff for each color
$c\in\Delta_d$, all $\hat c$-residues of $\Gamma$ represent closed $(d-1)$-manifolds (resp. the $(d-1)$-sphere).

In virtue of the bijection described in Remark \ref{correspondence-sing-boundary}, a $(d+1)$-colored  graph $(\Gamma,\gamma)$ is said to {\it represent}
a compact PL $d$-manifold $M$ with no spherical boundary components if and only if  it represents the associated singular manifold $\widehat M$.

\begin{definition} \label{crystallization}
A {\it crystallization} of a closed PL $d$-manifold $M^d$ is a $(d+1)$-colored graph representing $M^d$, such that each $d$-residue is connected (i.e.  $g_{\hat i}=1$ \ $\forall i \in \Delta_d$).
\end{definition}

The following theorem extends to singular manifolds a well-known result - due to Pezzana (\cite{[Pezzana]}) - founding the combinatorial representation theory for closed PL-manifolds of arbitrary dimension via colored graphs (the so called {\it crystallization theory}).  See also \cite{Cristofori-Fomynikh-Mulazzani-Tarkaev} and \cite{Cristofori-Mulazzani} for the 3-dimensional case.

\begin{theorem} {\em (\cite[Theorem 1]{Casali-Cristofori-Grasselli})} \label{Theorem_gem}
Any singular $d$-manifold  - or, equivalently, any compact $d$-manifold with no spherical boundary components -  admits a $(d+1)$-colored graph representing it.

In particular, each closed PL $d$-manifold admits a crystallization.
\end{theorem}

\smallskip 

It is well known the existence of a particular set of embeddings of a bipartite (resp. non-bipartite) $(d+1)$-colored graph into orientable (resp. non-orientable) surfaces.

\begin{theorem}{\em (\cite{Gagliardi 1981})}\label{reg_emb}
Let $(\Gamma,\gamma)$ be a bipartite (resp. non-bipartite) $(d+1)$-colored graph of order $2p$. Then for each cyclic permutation $\varepsilon = (\varepsilon_0,\ldots,\varepsilon_d)$ of $\Delta_d$, up to inverse, there exists a cellular embedding, called \emph{regular}, of $(\Gamma,\gamma)$ into an orientable (resp. non-orientable) closed surface $F_{\varepsilon}(\Gamma)$
whose regions are bounded by the images of the $\{\varepsilon_j,\varepsilon_{j+1}\}$-colored cycles, for each $j \in \mathbb Z_{d+1}$.
Moreover, the genus (resp. half the genus)  $\rho_{\varepsilon} (\Gamma)$ of $F_{\varepsilon}(\Gamma)$ satisfies

\begin{equation*}
\chi (F_\varepsilon(\Gamma)) =  2 - 2\rho_\varepsilon(\Gamma)= \sum_{j\in \mathbb{Z}_{d+1}} g_{\varepsilon_j\varepsilon_{j+1}} + (1-d)p.
\end{equation*}

No regular embeddings of $(\Gamma,\gamma)$ exist into non-orientable (resp. orientable) surfaces.
\end{theorem}

\smallskip 

The \emph{Gurau degree} (often called {\it degree} in the tensor models literature) and the {\it regular genus} of a colored graph are defined in terms of the embeddings of Theorem \ref{reg_emb}.

\begin{definition} \label{G-degree}
Let $(\Gamma,\gamma)$ be a $(d+1)$-colored graph.
If $\{\varepsilon^{(1)}, \varepsilon^{(2)}, \dots , \varepsilon^{(\frac {d!} 2)}\}$ is the set of all cyclic permutations of $\Delta_d$ (up to inverse), $ \rho_{\varepsilon^{(i)}}$ ($i=1, \dots , \frac {d!} 2$) is called the \emph{regular genus of $\Gamma$ with respect to the permutation $\varepsilon^{(i)}$}. Then, the \emph{Gurau degree} (or \emph{G-degree} for short) of $\Gamma$, denoted by  $\omega_{G}(\Gamma)$, is defined as
\begin{equation*}
 \omega_{G}(\Gamma) \ = \ \sum_{i=1}^{\frac {d!} 2} \rho_{\varepsilon^{(i)}}(\Gamma)
\end{equation*}
and the \emph{regular genus} of $\Gamma$, denoted by $\rho(\Gamma)$, is defined as
\begin{equation*}
 \rho(\Gamma) \ = \ \min\, \Big\{\rho_{\varepsilon^{(i)}}(\Gamma)\ /\ i=1,\ldots,\frac {d!} 2\Big\}.
\end{equation*}
\end{definition}

Note that, in dimension $2$, any bipartite (resp. non-bipartite) $3$-colored graph $(\Gamma,\gamma)$ represents an orientable (resp. non-orientable) surface $|K(\Gamma)|$ and $\rho(\Gamma)= \omega_G(\Gamma)$ is exactly the genus (resp. half the genus) of $|K(\Gamma)|.$
On the other hand, for $d\geq 3$, the G-degree of any $(d+1)$-colored graph (resp. the regular genus of any $(d+1)$-colored graph representing a closed PL $d$-manifold) is proved to be a non-negative {\it integer}, both in the bipartite and non-bipartite case:  see \cite[Proposition 7]{Casali-Cristofori-Dartois-Grasselli} (resp. \cite[Proposition A]{Chiavacci-Pareschi}).

\section{Proof of the general results} 
 \label{sec: Factorizations}

Within combinatorics,  the problem of the existence of $m$-cycle decompositions of the complete graph $K_n$, or of the complete multigraph $\lambda K_n$ (i.e. the multigraph with $n$ vertices and with $\lambda$ edges joining each pair of distinct vertices) is long standing: a survey result, for general $m$, $n$ and $\lambda$,  is given in \cite[Theorem 1.1]{[Bryant&al]}.

Moreover, the following results hold, concerning Hamiltonian cycles (i.e. $m=n$) in $K_n$, both in the case $n$ odd and in the case $n$ even.

\begin{proposition} \ {\rm \cite[Theorem 1.3]{[Bryant]}} \label{Bryant}
For all odd $n \ge 3$ there exists a partition of all Hamiltonian cycles
of $K_n$ into $(n-2)!$  Hamiltonian cycle decompositions of $K_n$.
\end{proposition}

\begin{proposition} \ {\rm \cite[Theorem 2.2]{Zhao-Kang}} \label{Zhao-Kang}
For all even $n \ge 4$ there exists a partition of all Hamiltonian cycles of $K_n$ into $\frac{(n-2)!}{2}$ classes, so that each edge of $K_n$ appears in exactly two cycles belonging to the same class.
\end{proposition}

Figure 1 describes - as an example of Proposition \ref{Bryant}  - the six Hamiltonian cycle decompositions of $K_5$, each of them containing a pair of disjoint Hamiltonian cycles (given by the dashed and continuous edges respectively).  
Note that, by labelling the vertices of $K_n$ with the elements of $\Delta_{n-1}$, each Hamiltonian cycle in $K_n$ defines a cyclic permutation of $\Delta_{n-1}$ (together with its inverse), and viceversa. 

\begin{figure} \label{figure1} 
\centerline{\scalebox{0.4}{\includegraphics{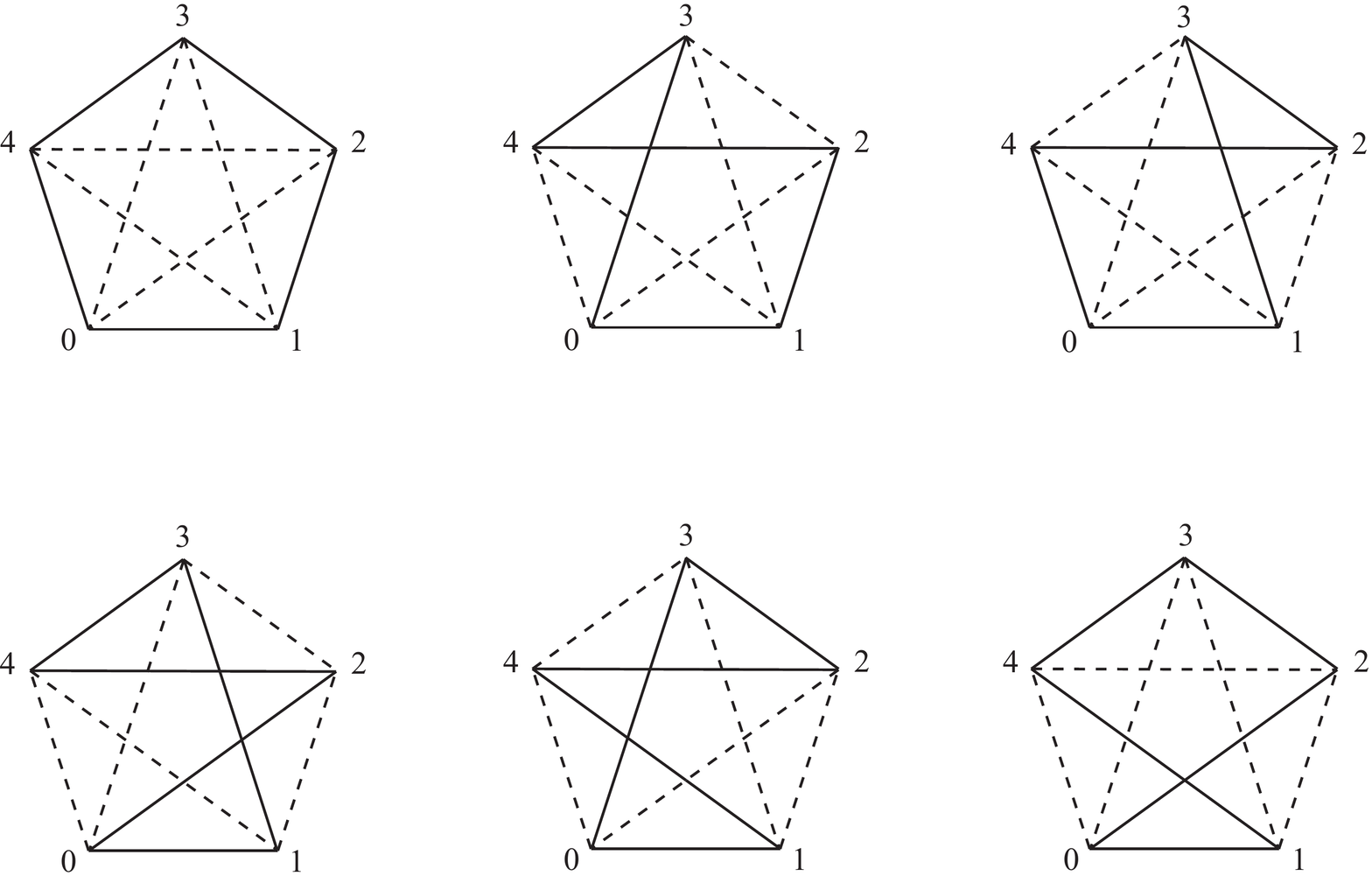}}}
\caption{The six Hamiltonian cycle decompositions of $K_5$}
\end{figure}

 \bigskip

On the other hand, the following statement regarding the G-degree has been recently proved.

\begin{proposition} \ {\rm \cite[Proposition 7]{Casali-Cristofori-Dartois-Grasselli}} \label{ridotto G-degree}
If $(\Gamma, \gamma)$ is a $(d+1)$-colored graph  of order $2p$ ($d\geq 3$), then
\begin{equation} \label{formula_G-degree}
\omega_G(\Gamma) \ = \ \frac{(d-1)!}{2} \cdot \Big(d + p \cdot (d-1) \cdot \frac d 2  -  \sum_{r,s \in \Delta_d} g_{rs}\Big).\end{equation}
As a consequence, the G-degree of any $(d+1)$-colored graph is a non-negative integer multiple of $\frac {(d-1)!} 2$.
\end{proposition}

The result of Proposition \ref{ridotto G-degree}, which was originally stated in the bipartite case (see \cite{BGR}), suggested the definition, for $d\geq 3$, of the (integer) {\it reduced G-degree}
$$\omega^{\prime}_G(\Gamma) \ = \ \frac {2} {(d-1)!} \cdot \omega_G(\Gamma),$$
which is used by many authors within tensor models theory (see for example \cite{Gurau-book}).\footnote{In fact, the exponents of $N^{-1}$ in the  $1/N$ expansion of formula \eqref{1/N expansion} are all (non-negative) integers $\omega_G^{\prime}$.}

Actually, we are able to prove that, if $d \ge 4$ is even, under 
rather weak hypotheses, 
the G-degree is multiple of (d-1)! (or, equivalently, the reduced G-degree is even).

\begin{proposition} \label{ridottissimo G-degree}
If $d \ge 4$ is even, and  $(\Gamma, \gamma)$ is a $(d+1)$-colored graph such that 
each $3$-residue is bipartite and each $d$-residue is either bipartite or non-bipartite with integer regular genus with respect to any permutation, 
then
$$\omega_G(\Gamma) \equiv 0 \mod (d-1)!$$
\end{proposition}

\dimo 
Let $(\Gamma,\gamma)$ be a $(d+1)$-colored graph, with $d \ge 4,$ $d$ even.
Since $n = d+1$ is odd, Proposition \ref{Bryant} implies that all $\frac{d!}{2} $ cyclic permutations (up to inverse) of $\Delta_d$ can be partitioned in $(d-1)!$ classes, each containing $d/2$ cyclic permutations,
$\bar \varepsilon^{(1)}, \bar \varepsilon^{(2)}, \dots , \bar \varepsilon^{(\frac {d} 2)}$ say,
so that
$$ \sum_{i=1}^{d/2} \Big[ \sum_{j\in \mathbb{Z}_{d+1}} g_{\bar \varepsilon^{(i)}_j,\bar \varepsilon^{(i)}_{j+1}} \Big] = \sum_{r,s\in \Delta_{d}} g_{r,s}.$$
Hence, by Theorem \ref{reg_emb},
$$2 \cdot \frac d 2 - 2 \cdot \sum_{i=1}^{d/2} \rho_{\bar \varepsilon^{(i)}} =  \sum_{r,s \in \Delta_d} g_{rs}   + p \cdot (1-d) \cdot \frac d 2,$$
and
\begin{equation} \label{eq.sum}
2 \cdot  \sum_{i=1}^{d/2} \rho_{\bar \varepsilon^{(i)}} =  d + p \cdot (d-1) \cdot \frac d 2  -  \sum_{r,s \in \Delta_d} g_{rs}.
\end{equation}
This proves that the quantity  $\sum_{i=1}^{d/2} \rho_{\bar \varepsilon^{(i)}}$ is constant, for each class of the above partition; then,
$$\omega_G (\Gamma)=  (d-1)! \cdot \sum_{i=1}^{d/2} \rho_{\bar \varepsilon^{(i)}}$$
immediately follows.\footnote{In this way, one can reobtain - for $d \ge 4$ even - relation \eqref{formula_G-degree}.} 

Now, the hypotheses on the $3$-residues and $d$-residues 
of $\Gamma$ directly implies, in virtue of  Theorem \ref{reg_emb} and  of \cite[Lemma 2]{Casali-Cristofori-Grasselli} (originally proved in \cite[Lemma 4.2]{Chiavacci-Pareschi}),
that  $\sum_{i=1}^{d/2} \rho_{\bar \varepsilon^{(i)}}$ is an integer, and hence $\omega_G (\Gamma) \equiv 0 \mod (d-1)!$ \ \ (or equivalently, $\omega_G^{\prime} (\Gamma) = 2 \cdot \frac{\omega_G (\Gamma)}{(d-1)!}$ is even).
\qed

\medskip

Then, the results stated in Section \ref{intro} trivially follow. 

\medskip

\noindent {\it Proof  of Theorems 1 and 2.}
\par \noindent
Both in the case of $(\Gamma, \gamma)$ bipartite and in the case of $(\Gamma, \gamma)$ representing a singular $d$-manifold, 
the residues of $\Gamma$ obviously satisfy the hypotheses of Proposition \ref{ridottissimo G-degree}.   
Hence, if $d\ge 4$ is even,  $\omega_G(\Gamma) \equiv 0 \mod (d-1)!$ holds. \
\qed

\bigskip

Another particular situation is covered by Proposition  \ref{ridottissimo G-degree}, as the following corollary explains.

\begin{corollary}
Let $(\Gamma,\gamma)$ be a $(d+1)$-colored graph, with $d \ge 4,$ $d$ even.
If $(\Gamma,\gamma)$  is a non-bipartite $(d+1)$-colored graph such that each  $d$-residue is bipartite,
then $$\omega_G (\Gamma) \equiv 0 \mod (d-1)!$$
\vskip-0.3truecm \ \ \qed
\end{corollary}

Note that there exist $(d+1)$-colored graphs, with $d$ even, $d \ge 4$ and odd reduced G-degree: of course, in virtue of Theorems 1 and 2, the represented $d$-pseudomanifold must be non-orientable and it can't be a singular $d$-manifold. As an example, for each $d \ge 4$, the $(d+1)$-colored graph  $(\Gamma,\gamma)$ of Figure 2 represents the $(d-2)$-th suspension $\Sigma$ of the real projective plane $\mathbb{RP}^2$  and $\omega_G^{\prime}(\Gamma)=d-1.$ 
It is easy to check that $\Gamma$ is non-bipartite and $\Sigma$ is not a singular manifold, since all $d$-residues of $\Gamma$, with the exception of $\Gamma_{\hat 0}$, do not represent a closed $(d-1)$-manifold. 
Moreover, its (unique) $\{0,1,c \}$-residue, for any $c \in \{2, \dots, d\}$, is non-bipartite and represents  the  non-orientable genus one surface $\mathbb{RP}^2$.

\begin{figure} \label{figure2} 
\centerline{\scalebox{0.4}{\includegraphics{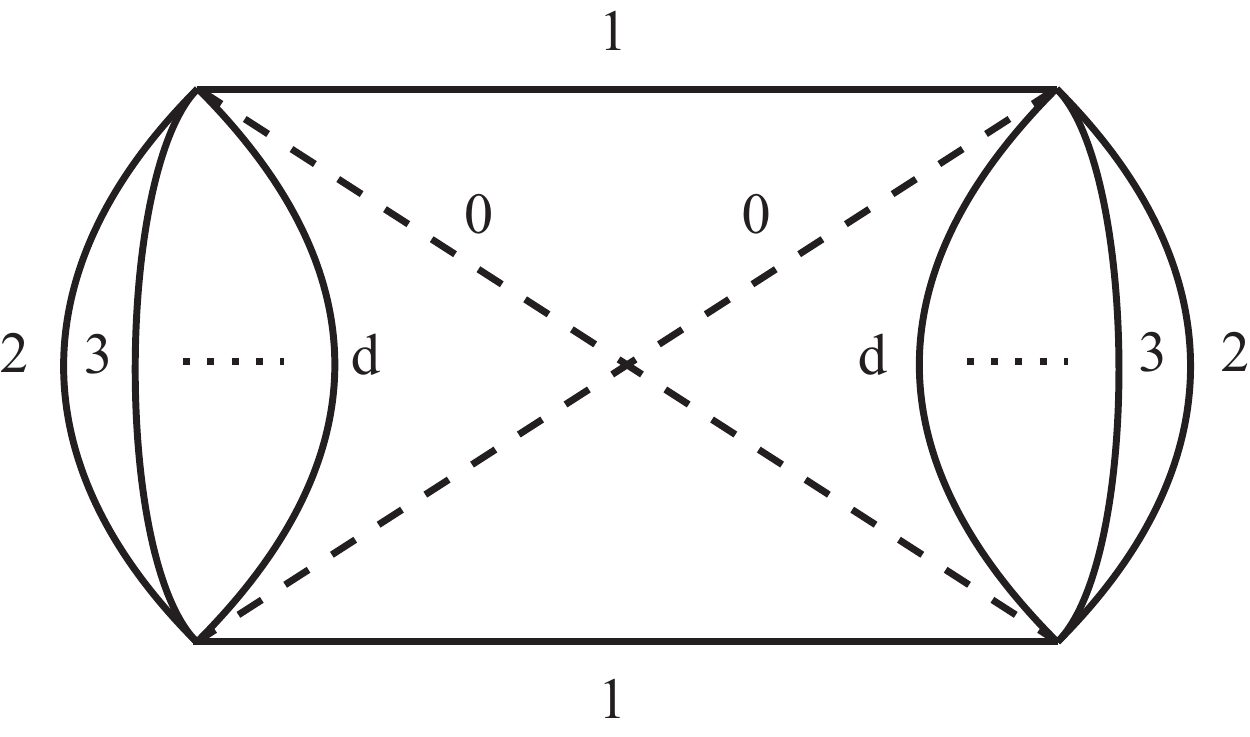}}}
\caption{A $(d+1)$-colored graph representing the $(d-2)$-th suspension of 
$\mathbb{RP}^2$}
\end{figure}

\begin{remark} \label{rem:odd} {\rm
In the case $d \ge 3$ odd, Proposition \ref{Zhao-Kang} implies that all $\frac{d!}{2} $ cyclic permutations (up to inverse) of $\Delta_d$ can be partitioned in $\frac{(d-1)!}{2}$ classes, each containing $d$ cyclic permutations,
$\bar \varepsilon^{(1)}, \bar \varepsilon^{(2)}, \dots , \bar \varepsilon^{(d)}$ say,
so that
$$ \sum_{i=1}^{d} \Big[ \sum_{j\in \mathbb{Z}_{d+1}} g_{\bar \varepsilon^{(i)}_j,\bar \varepsilon^{(i)}_{j+1}} \Big] = 2 \sum_{r,s\in \Delta_{d}} g_{r,s}.$$
Hence, a reasoning similar to the one used to prove Proposition \ref{ridottissimo G-degree} yields an alternative proof - for $d \ge 3$ odd - of relation \eqref{formula_G-degree}: 
since $$2 \cdot d - 2 \cdot \sum_{i=1}^d \rho_{\bar \varepsilon^{(i)}} = 2 \sum_{r,s \in \Delta_d} g_{rs}   + p \cdot (1-d) \cdot d$$
and
\begin{equation} \label{eq.sum(odd)}  \sum_{i=1}^d \rho_{\bar \varepsilon^{(i)}}  =  d + p \cdot (d-1) \cdot \frac d 2  -  \sum_{r,s \in \Delta_d} g_{rs}
\end{equation}
hold, then
$$\omega_G (\Gamma)= \frac{(d-1)!}{2} \cdot \sum_{i=1}^d \rho_{\bar \varepsilon^{(i)}} = \frac{(d-1)!}{2} \cdot \Big(d + p \cdot (d-1) \cdot  \frac d 2 -  \sum_{r,s \in \Delta_d} g_{rs} \Big)$$
directly follows.}
\end{remark}

\begin{remark} \label{rem:constant-general} {\rm
It is worthwhile to stress that, for $d$ even (resp. odd), formula \eqref{eq.sum} of Proposition \ref{ridottissimo G-degree} (resp. formula \eqref{eq.sum(odd)} of Remark \ref{rem:odd}) proves that the sum $\sum_{i=1}^{d/2} \rho_{\bar \varepsilon^{(i)}}$ (resp. $ \sum_{i=1}^d \rho_{\bar \varepsilon^{(i)}}$) of all regular genera with respect to the $d/2$ (resp. $d$) permutations belonging to the same class
is half the reduced G-degree (resp. is the reduced G-degree)
$$ \omega^{\prime}_G(\Gamma) \ = \ d + p \cdot (d-1) \cdot \frac d 2  -  \sum_{r,s \in \Delta_d} g_{rs}, $$  
i.e. it is a constant which does not depend on the chosen partition class.

Hence, the regular genus $\rho(\Gamma)$ of the graph $\Gamma$ is realized by the (not necessarily unique) permutation  $\varepsilon$ which maximizes the difference
$$ \rho_{\hat \varepsilon}(\Gamma) - \rho_{\varepsilon}(\Gamma),$$
where $ \rho_{\hat \varepsilon}(\Gamma)$ denotes the sum of the genera with respect to all other permutations of the same partition class.
}
\end{remark}

\section{The $4$-dimensional case}  \label{sec: n=4}

In the $4$-dimensional setting, the above combinatorial properties allow to prove further results about the G-degree.

In fact, it is easy to check that, for each  cyclic permutation $\varepsilon = (\varepsilon_0, \varepsilon_1, \varepsilon_2, \varepsilon_3, \varepsilon_4)$ of $\Delta_4$, we have
\begin{equation} \label{2-factorization n=4}
\{ (\varepsilon_j,\varepsilon_{j+1}) \ / \ j\in \mathbb{Z}_{5}\} \  \cup  \ \{ (\varepsilon_j, \varepsilon_{j+2}) \ / j\in \mathbb{Z}_{5}\} \ = \  \{ (i,j) \ / \  i,j\in \Delta_{5}, i \ne j\}.
\end{equation}

Let us denote by $\varepsilon^{\prime}=(\varepsilon^{\prime}_0, \varepsilon^{\prime}_1, \varepsilon^{\prime}_2, \varepsilon^{\prime}_3, \varepsilon^{\prime}_4)$  the permutation $(\varepsilon_0, \varepsilon_2, \varepsilon_4, \varepsilon_1, \varepsilon_3)$ of $\Delta_4$, which may be said to be {\it associated} to $\varepsilon$.

Note that, when $d=4$, the only partition of all $12$ cyclic permutations of $\Delta_4$ (up to inverse) is given by the six classes containing a given permutation $\varepsilon$ and its associated $\varepsilon^{\prime}$: see Figure 1, where each class of the Hamiltonian cycle decomposition of $K_5$ ensured by Proposition \ref{Bryant}  is shown to correspond to a pair $(\varepsilon_i, \varepsilon_i^{\prime})$ of associated permutations, for $i \in \mathbb N_6.$ For example, the first class corresponds to the identical permutation $\varepsilon_1 = (0,1,2,3,4)$ and its associated permutation $\varepsilon_1^{\prime} = (0,2,4,1,3).$ 

\smallskip

Then, the following result holds.

\begin{proposition} \label{iper-reduced G-degree(d=4)}
For each $5$-colored graph  $(\Gamma, \gamma),$  and for each pair ($\varepsilon, \varepsilon^{\prime})$ of associated cyclic permutations of $\Delta_4$,
$$ \omega_{G}(\Gamma)  \ = \ 6 \Big(\rho_{\varepsilon}(\Gamma) + \rho_{\varepsilon^{\prime}}(\Gamma)\Big).$$
 \end{proposition}

\dimo 
Equality \eqref{2-factorization n=4} directly yields
$$ \sum_{j\in \mathbb{Z}_{5}} g_{\varepsilon_j,\varepsilon_{j+1}} +  \sum_{j\in \mathbb{Z}_{5}} g_{\varepsilon^{\prime}_j,\varepsilon^{\prime}_{j+1}}= \sum_{i,j\in \Delta_5} g_{i,j}.$$

As a consequence, the sum of all regular genera of $\Gamma$ with respect to the $12$ cyclic permutations (up to inverse) of $\Delta_4$  is six times the sum between the regular genera of $\Gamma$ with respect to any pair $\varepsilon$, $\varepsilon^{\prime}$ of associated permutations:
$$ \omega_{G}(\Gamma)  \ = \ 6 \Big(\rho_{\varepsilon}(\Gamma) + \rho_{\varepsilon^{\prime}}(\Gamma)\Big).$$
\vskip-0.6truecm \ \ \qed

\begin{remark} \label{sum_constant} {\rm
By Proposition \ref{iper-reduced G-degree(d=4)}, for any 5-colored graph the sum between the regular genera of $\Gamma$ with respect to any pair $\varepsilon$, $\varepsilon^{\prime}$ of associated cyclic permutations  is
constant (see equality \eqref{eq.sum} and Remark \ref{rem:constant-general}, for $d=4$):
$$ \rho_{\varepsilon}(\Gamma) + \rho_{\varepsilon^{\prime}}(\Gamma) \ = \ \frac 1 2  \omega'_G(\Gamma) =  2+3p - \frac 1 2  \sum_{i,j\in \Delta_{5}} g_{i,j}.$$
Hence,
the regular genus $\rho(\Gamma)$ of the graph $\Gamma$ is realized by the (not necessarily unique) permutation $\varepsilon$ so that $\rho_{\varepsilon^{\prime}}(\Gamma) - \rho_{\varepsilon}(\Gamma)$ is maximal.
}
\end{remark}

Moreover:
\begin{proposition} \label{differences}
\begin{itemize}
\item[(a)]
If  $(\Gamma, \gamma)$ is a $5$-colored graph, then for each pair ($\varepsilon, \varepsilon^{\prime})$ of associated cyclic permutations of $\Delta_4$,
$$ 2 \Big(\rho_{\varepsilon^{\prime}}(\Gamma) - \rho_{\varepsilon}(\Gamma) \Big)  \ = \   \sum_{j\in \mathbb{Z}_{5}} g_{\varepsilon_j,\varepsilon_{j+1}} - \sum_{j\in \mathbb{Z}_{5}} g_{\varepsilon_j,\varepsilon_{j+2}}$$
\item[(b)]
If  $(\Gamma, \gamma)$ is a $5$-colored graph representing a singular $4$-manifold $M^4$,  then for each pair $(\varepsilon, \varepsilon^{\prime})$ of associated cyclic permutations of $\Delta_4$,
$$ \rho_{\varepsilon^{\prime}}(\Gamma) - \rho_{\varepsilon}(\Gamma)  \ = \   \sum_{j\in \mathbb{Z}_{5}} g_{\varepsilon_j,\varepsilon_{j+1},\varepsilon_{j+2}} - \sum_{j\in \mathbb{Z}_{5}} g_{\varepsilon_j,\varepsilon_{j+2}, \varepsilon_{j+4}}$$
\end{itemize}
\end{proposition}

\dimo 
Statement (a) is an easy consequence of Theorem \ref{reg_emb}:
$$2-2 \rho_{\varepsilon}(\Gamma) = \sum_{j\in \mathbb{Z}_{5}} g_{\varepsilon_j,\varepsilon_{j+1}} -3p$$
\noindent and
$$ 2-2 \rho_{\varepsilon^{\prime}}(\Gamma) = \sum_{j\in \mathbb{Z}_{5}} g_{\varepsilon^{\prime}_j,\varepsilon^{\prime}_{j+1}} -3p=  \sum_{j\in \mathbb{Z}_{5}} g_{\varepsilon_j,\varepsilon_{j+2}} -3p.$$   

On the other hand, relation $2  g_{r,s,t} = g_{r,s} + g_{r,t}  + g_{s,t} - p$
is known to be true for each order $2p$ $5$-colored graph representing a singular $4$-manifold (see \cite[Lemma 21]{Casali-Cristofori-Dartois-Grasselli}. As a consequence we have:
$$ 2 \sum_{j\in \mathbb{Z}_{5}} g_{\varepsilon_j,\varepsilon_{j+1},\varepsilon_{j+2}}= \sum_{r,s\in \mathbb{Z}_{5}} g_{r,s} + \sum_{j\in \mathbb{Z}_{5}} g_{\varepsilon_j,\varepsilon_{j+1}} - 5p$$  and
$$ 2 \sum_{j\in \mathbb{Z}_{5}} g_{\varepsilon_j,\varepsilon_{j+2},\varepsilon_{j+4}}= \sum_{r,s\in \mathbb{Z}_{5}} g_{r,s} + \sum_{j\in \mathbb{Z}_{5}} g_{\varepsilon_j,\varepsilon_{j+2}} - 5p.$$
By making the difference,
$$2 \Big(\sum_{j\in \mathbb{Z}_{5}} g_{\varepsilon_j,\varepsilon_{j+1},\varepsilon_{j+2}} - \sum_{j\in \mathbb{Z}_{5}} g_{\varepsilon_j,\varepsilon_{j+2},\varepsilon_{j+4}}\Big) =  \sum_{j\in \mathbb{Z}_{5}} g_{\varepsilon_j,\varepsilon_{j+1}} - \sum_{j\in \mathbb{Z}_{5}} g_{\varepsilon_j,\varepsilon_{j+2}}$$ is obtained; so, statement (b) follows, via statement (a).
\qed

Proposition \ref{differences} enables to obtain
the following improvement of  \cite[Proposition 29(a)]{Casali-Cristofori-Dartois-Grasselli}.

\begin{corollary}
Let $(\Gamma,\gamma)$ be a $5$-colored graph. \ Then:
$$ \omega_G(\Gamma) \ = \ 12 \cdot  \rho (\Gamma) \ \ \   \Longleftrightarrow$$
$$ \sum_{j\in\mathbb{Z}_{5}} g_{\varepsilon_j, \varepsilon_{j+1}} =  \sum_{j\in\mathbb{Z}_{5}} g_{\varepsilon_j, \varepsilon_{j+2}} \ \text{for each cyclic permutation} \ \varepsilon \ \text{of}  \ \Delta_4.$$
\end{corollary}

\dimo 
From \cite[Proposition 29 (a)]{Casali-Cristofori-Dartois-Grasselli}, it is known that
$$ \omega_G(\Gamma) \ = \ 12 \cdot  \rho (\Gamma) \ \ \   \Longrightarrow \ \ \ \  \sum_{j \in \mathbb{Z}_{5}} g_{\bar \varepsilon_j, \bar \varepsilon_{j+1}} =  \sum_{j\in\mathbb{Z}_{5}} g_{\bar \varepsilon_j, \bar \varepsilon_{j+2}},$$
$\bar \varepsilon$ being the cyclic permutation of $\Delta_4$ such that $\rho(\Gamma) = \rho_{\bar \varepsilon}(\Gamma).$
So, $\rho_{\bar \varepsilon^{\prime}}(\Gamma) - \rho_{\bar \varepsilon}(\Gamma) = 0$ directly follows via Proposition \ref{differences} (a).
Now, since $\rho_{\varepsilon}(\Gamma) + \rho_{\varepsilon^{\prime}}(\Gamma)$ is constant for each pair $(\varepsilon, \varepsilon^{\prime})$ of associated cyclic permutations of $\Delta_4$ (see Remark \ref{sum_constant}), then:
$$|\rho_{\varepsilon^{\prime}}(\Gamma) - \rho_{\varepsilon}(\Gamma)| \le \rho_{\bar \varepsilon^{\prime}}(\Gamma) - \rho_{\bar \varepsilon}(\Gamma).$$
Hence, $\rho_{\bar \varepsilon^{\prime}}(\Gamma) - \rho_{\bar \varepsilon}(\Gamma) = 0$ implies $\sum_{j\in\mathbb{Z}_{5}} g_{\varepsilon_j, \varepsilon_{j+1}} =  \sum_{j\in\mathbb{Z}_{5}} g_{\varepsilon_j, \varepsilon_{j+2}}$ (i.e. $\rho_{\varepsilon^{\prime}}(\Gamma) - \rho_{\varepsilon}(\Gamma) = 0$, in virtue of Proposition \ref{differences} (a)) for each cyclic permutation $\varepsilon$ of $\Delta_4$.

\par \noindent
The reversed implication is straightforward, via Proposition \ref{differences}(a).
\qed

\begin{proposition} \label{inequalities Euler characteristic}
If $(\Gamma, \gamma)$  is an order $2p$ 5-colored graph representing
a singular 4-manifold $M^4$, then, for each pair ($\varepsilon, \varepsilon^{\prime})$ of associated cyclic permutations of $\Delta_4$:
$$ \chi(M^4)  \ = \ \Big(\rho_{\varepsilon}(\Gamma) + \rho_{\varepsilon^{\prime}}(\Gamma)\Big) -p + \sum_{i\in \Delta_{4}}  g_{\hat i} -2.$$
 \end{proposition}

\dimo 
It is sufficient to apply  Proposition \ref{iper-reduced G-degree(d=4)}  to the third equality of \cite[Proposition 22]{Casali-Cristofori-Dartois-Grasselli}.
\qed

\smallskip 

Let us now recall two particular types of crystallizations introduced and studied in \cite{[Basak-Casali 2016]} and \cite{[Basak]}\footnote{Both semi-simple and weak semi-simple crystallizations generalize the notion of {\it simple crystallizations} for simply-connected PL $4$-manifolds: see \cite{[Basak-Spreer 2016]} and \cite{[Casali-Cristofori-Gagliardi JKTR 2015]}.}: they are proved to be ``minimal" with respect to regular genus, among all graphs representing the same PL $4$-manifold.

\begin{definition} {\rm A crystallization of a PL $4$-manifold $M^4$ with \, $rk(\pi_1(M^4))= m \ge 0$ is called a {\em semi-simple crystallization}
if \ $g_{j,k,l} = 1 + m \ \ \forall \ j,k,l \in \Delta_4.$

A crystallization of a PL $4$-manifold $M^4$ with \, $rk(\pi_1(M^4))= m $ is called a {\em weak semi-simple crystallization}
if \ $g_{i, i+1, i +2} = 1 + m \ \ \forall \ i \in \Delta_4$ (where the additions in subscripts are intended in $\mathbb{Z}_{5}$).
}
\end{definition}

According to \cite{[Basak-Casali 2016]}, for each order $2p$ crystallization  $(\Gamma, \gamma)$ of a closed  PL 4-manifold $M^4$, with \, $rk(\pi_1(M^4))= m0$, let us set
$$g_{j,k,l}= 1+m + t_{j,k,l},  \ \ \text{with} \ t_{j,k,l}  \ \ \forall j,k,l \in \Delta_4.$$
Semi-simple (resp. weak semi-simple) crystallizations turn out to be characterized by $t_{j,k,l} = 0$  \ $\forall j,k,l \in \Delta_4$ (resp. $t_{i, i+1, i +2} = 0$ \ $\forall i \in \mathbb{Z}_{5}$).

In \cite{[Basak-Casali 2016]}, the relation
\begin{equation} \label{p and q}
p= 3 \chi(M^4)+ 5(2m -1) +  \sum_{j,k,l \in \Delta_4} t_{j,k,l}
\end{equation}
is proved to hold; hence, $p = \bar p + q$ follows, where $q= \sum_{j,k,l \in \Delta_4} t_{j,k,l} \ge 0$ and $\bar p= 3 \chi(M^4)+ 5(2m -1)$ is the minimum possible half order of a crystallization of $M^4$, which is attained if and only if $M^4$ admits semi-simple crystallizations.

With the above notations, the following results can be obtained.

\begin{proposition} \label{sum/difference}
Let $(\Gamma,\gamma)$ be an order $2p$ crystallization  of a closed PL 4-manifold $M^4$, with $rk(\pi_1(M^4))=m.$ Then, for each pair ($\varepsilon, \varepsilon^{\prime})$ of associated cyclic permutations of $\Delta_4$, with $\rho_{\varepsilon}(\Gamma) \le \rho_{\varepsilon^{\prime}}(\Gamma)$:
\begin{equation}
\label{eq.difference}
\rho_{\varepsilon^{\prime}}(\Gamma) - \rho_{\varepsilon}(\Gamma) = q - 2 \sum_{i\in \mathbb{Z}_{5}}  t_{\varepsilon_i,\varepsilon_{i+2},\varepsilon_{i+4}} \ \le \ q.
\end{equation}
Moreover, The following statements are equivalent:
\begin{itemize}
\item[(a)] a cyclic permutation $\varepsilon$ of $\Delta_4$ exists, such that $\rho_{\varepsilon^{\prime}}(\Gamma) - \rho_{\varepsilon}(\Gamma) =q;$
\item[(b)]  $\Gamma$ is a weak semi-simple crystallization;
\item[(c)] $ \rho(\Gamma)= 2\chi(M^4) + 5 m -4$.\footnote{Note that, in this case, the permutation $\bar \varepsilon$ such that $\rho(\Gamma)= \rho_{\bar \varepsilon}(\Gamma)$ coincides with the permutation $\varepsilon$ of point (a); moreover, $\rho(\Gamma)= \mathcal G(M^4)$ holds, in virtue of the inequality  $\mathcal G(M^4) \ge 2\chi(M^4) + 5 m -4$, proved in  \cite{[Basak-Casali 2016]}.}
\end{itemize}
\end{proposition}

\dimo 
In virtue of Proposition \ref{differences}(b), an easy computation proves relation \eqref{eq.difference}:
$$ \begin{aligned} \rho_{\varepsilon^{\prime}}(\Gamma) - \rho_{\varepsilon}(\Gamma) \ & = \  \sum_{i\in \mathbb{Z}_{5}}  (m+1 + t_{\varepsilon_i,\varepsilon_{i+1},\varepsilon_{i+2}} ) - \sum_{i\in \mathbb{Z}_{5}} (m+1 + t_{\varepsilon_i,\varepsilon_{i+2},\varepsilon_{i+4}}) \ = \\
\ & = \   \sum_{i\in \mathbb{Z}_{5}}  t_{\varepsilon_i,\varepsilon_{i+1},\varepsilon_{i+2}} - \sum_{i\in \mathbb{Z}_{5}}   t_{\varepsilon_i,\varepsilon_{i+2},\varepsilon_{i+4}}  \ = \\
\ & = \  \sum_{j,k,l \in \Delta_4}  t_{j,k,l} - 2 \sum_{i\in \mathbb{Z}_{5}}  t_{\varepsilon_i,\varepsilon_{i+2},\varepsilon_{i+4}} \ = \\
\ & = \  q - 2 \sum_{i\in \mathbb{Z}_{5}}  t_{\varepsilon_i,\varepsilon_{i+2},\varepsilon_{i+4}} \ \le \  q. \end{aligned}$$

Now, if $\Gamma$ is a weak semi-simple crystallization, by definition itself a cyclic permutation $\varepsilon$ of $\Delta_4$ exists\footnote{$\varepsilon$ turns out to be the permutation of $\Delta_4$ associated to $\varepsilon^{\prime}=(0,1,2,3,4)$.}, so that, for each $i \in \mathbb Z_5,$ \ $g_{\varepsilon_i,\varepsilon_{i+2},\varepsilon_{i+4}}= m+1,$ i.e. $t_{\varepsilon_i,\varepsilon_{i+2},\varepsilon_{i+4}}=0.$  So, $\rho_{\varepsilon^{\prime}}(\Gamma) - \rho_{\varepsilon}(\Gamma) =q$ easily follows from relation \eqref{eq.difference}.

On the other hand,  if  a cyclic permutation $\varepsilon$ of $\Delta_4$ exists, so that $\rho_{\varepsilon^{\prime}}(\Gamma) - \rho_{\varepsilon}(\Gamma) =q,$  relation \eqref{eq.difference} yields  $t_{\varepsilon_i,\varepsilon_{i+2},\varepsilon_{i+4}}=0,$ i.e. $g_{\varepsilon_i,\varepsilon_{i+2},\varepsilon_{i+4}}= m+1,$ which is exactly - up to color permutation - the definition of weak semi-simple crystallization.
Hence, (a) and (b) are proved to be equivalent.  

Then, by comparing Proposition \ref{inequalities Euler characteristic} (with the assumption $g_{\hat i}=1$ \  $\forall i \in \Delta_4$) and relation \eqref{p and q}, we obtain:
\begin{equation} \label{sum_q}
\rho_{\varepsilon^{\prime}}(\Gamma) + \rho_{\varepsilon}(\Gamma) = 2 ( 2 \chi(M^4) + 5m -4) + q.
\end{equation}
By making use of relation \eqref{eq.difference},
\begin{equation} \label{rho-corrected}
 \rho_{\varepsilon}(\Gamma)  = 2 \chi(M^4) + 5m -4 + \sum_{i\in \mathbb{Z}_{5}}  t_{\varepsilon_i,\varepsilon_{i+2},\varepsilon_{i+4}}
\end{equation}
easily follows, as well as
\begin{equation} \label{rho-epsilon'} \rho_{\varepsilon^{\prime}}(\Gamma)  = 2 \chi(M^4) + 5m -4 + \sum_{i\in \mathbb{Z}_{5}}  t_{\varepsilon_i,\varepsilon_{i+1},\varepsilon_{i+2}}.
\end{equation}
Relation \eqref{rho-corrected} directly yields the co-implication between statements (b) and (c).
\qed

We conclude the paper with a list of remarks, which arise from the previous results.

\begin{remark} {\rm With the above notations, formulae \eqref{rho-corrected} and \eqref{rho-epsilon'} immediately give
$$ 2 + \frac {\rho_{\varepsilon}(\Gamma)}  2 - \frac {5 m} 2 - \frac q 4 \ \le \ \chi(M^4)  \ \le \ 2 + \frac {\rho_{\varepsilon^{\prime}}(\Gamma)} 2 - \frac {5 m} 2 - \frac q 4$$
for each crystallization  of a closed PL 4-manifold $M^4$.
Another double inequality concerning the Euler characteristic of a closed  PL 4-manifold $M^4$ and the regular genera of any order $2p$ crystallization of $M^4$ with respect to a pair of associated permutations may be easily obtained from  Proposition \ref{inequalities Euler characteristic},
by making use of the assumptions $\sum_{i\in \Delta_{4}}  g_{\hat i}=5$ and $\rho_{\varepsilon}(\Gamma) \le \rho_{\varepsilon^{\prime}}(\Gamma):$
$$ 2\rho_{\varepsilon}(\Gamma) -p +3 \ \le \ \chi(M^4)  \ \le \ 2 \rho_{\varepsilon^{\prime}}(\Gamma) -p +3.$$
Note that such double inequalities assume a specific relevance in case of ``low" difference between $\rho_{\varepsilon^{\prime}}(\Gamma)$ and $\rho_{\varepsilon}(\Gamma)$ (in particular if  $\rho_{\varepsilon^{\prime}}(\Gamma)= \rho_{\varepsilon}(\Gamma)$ occurs, possibly with $\rho(\Gamma) < \rho_{\varepsilon}(\Gamma)$).

On the other hand, Proposition \ref{inequalities Euler characteristic} and formula \eqref{eq.difference}  (resp. formula \eqref{sum_q}) give
$$\chi(M^4)  \ = \ 2 \rho_{\varepsilon}(\Gamma) -p + 3 + (q -2 \sum_{i\in \mathbb{Z}_{5}}  t_{\varepsilon_i,\varepsilon_{i+2},\varepsilon_{i+4}})$$
$$ \text{\Big(resp. \ \ }  \chi(M^4)  \ = \ 2 + \frac {\rho_{\varepsilon}(\Gamma)}  2 - \frac {5 m} 2 - \frac 1 2 \sum_{i\in \mathbb{Z}_{5}}   t_{\varepsilon_i,\varepsilon_{i+2},\varepsilon_{i+4}} \ \text{\Big)}.$$
Hence, the following double inequalities arise, too, both involving the regular genus with respect to only one cyclic permutation: $$ 2\rho_{\varepsilon}(\Gamma) -p +3 \ \le \ \chi(M^4)  \ \le \ 2 \rho_{\varepsilon}(\Gamma) - p + q +3; $$
$$ 2 + \frac {\rho_{\varepsilon}(\Gamma)}  2 - \frac {5 m} 2 - \frac q 4  \ \le \ \chi(M^4)  \ \le \ 2 + \frac {\rho_{\varepsilon}(\Gamma)}  2 - \frac {5 m} 2 .$$
}
\end{remark}

\begin{remark} {\rm Note that relation \eqref{sum_q} exactly corresponds, via Proposition \ref{iper-reduced G-degree(d=4)}, to \cite[Proposition 27]{Casali-Cristofori-Dartois-Grasselli}. Moreover, relation \eqref{rho-corrected}\footnote{Actually, relation \eqref{rho-corrected} corrects a trivial error in the proof (and statement) of \cite[Lemma 7]{[Basak]}, not affecting the implications in order to prove the main result of that paper.}
directly implies the inequality $\rho_{\varepsilon}(\Gamma)  \ge 2 \chi(M^4) + 5m -4$ (which is one of the upper bounds obtained in \cite{[Basak-Casali 2016]}) and ensures that, for each crystallization $(\Gamma,\gamma)$ of a PL 4-manifold, the regular genus $\rho(\Gamma)$ is realized by the (not necessarily unique) permutation $\varepsilon$ so that $\sum_{i\in \mathbb{Z}_{5}}  t_{\varepsilon_i,\varepsilon_{i+2},\varepsilon_{i+4}}$ (or, equivalently, $\sum_{i\in \mathbb{Z}_{5}}  g_{\varepsilon_i,\varepsilon_{i+2},\varepsilon_{i+4}}$) is minimal.
Finally, the equivalence between items (b) and (c) in the above Proposition \ref{sum/difference} gives a direct proof of \cite[Theorem 2]{[Basak]}.
}
\end{remark}

\begin{remark} \label{rem:intersection-classes} {\rm Semi-simple crystallizations turn out to be the intersection (characterized by $q=0$) between the two classes of  weak semi-simple crystallizations and of crystallizations satisfying $\omega_G(\Gamma) \ = \ 12 \cdot  \rho (\Gamma)$ (and hence $\rho_{\varepsilon}(\Gamma) = 2 \chi(M^4) + 5m -4 + \frac 1 2 q,$ for each cyclic permutation $\varepsilon$ of $\Delta_4$).

Moreover, it is easy to check that $$ q\le 2  \ \ \ \ \Longrightarrow \ \ \ \ \Gamma \ \text{ is a weak semi-simple crystallization}.$$
In fact, if $q = \sum_{j,k,l \in \Delta_4} t_{j,k,l} \le 2, $ at most two triads $(j,k,l)$ of distinct elements in $\Delta_4$ exist, so that $g_{j,k,l}= 1+m + t_{j,k,l} > 1+m.$ This ensures the existence of a cyclic permutation $\varepsilon$ of $\Delta_4$ so that, for each $i \in \mathbb Z_5,$ \ $g_{\varepsilon_i,\varepsilon_{i+1},\varepsilon_{i+2}}= m+1,$ which is exactly the requirement for a weak semi-simple crystallization.
}
\end{remark}

\begin{remark} {\rm The formula obtained in \cite[Lemma 4.2]{Gurau-Ryan} for bipartite $(d+1)$-colored graphs and extended to the general case in \cite[Lemma 13]{Casali-Cristofori-Dartois-Grasselli} gives, if $d=4$,
$$\omega_G(\Gamma)= 3 \Big(p + 4 - \sum_{i\in \Delta_{4}}  g_{\hat i}\Big) + \sum_{i\in \Delta_{4}} \omega_G(\Gamma_{\hat i}),$$
where, for each  $i\in \Delta_{4}$, $\omega_G(\Gamma_{\hat i})$ denotes the sum of the G-degrees of the connected components of $\Gamma_{\hat i}.$
Hence, $\sum_{i\in \Delta_{4}} \omega_G(\Gamma_{\hat i})$ is always a multiple of $3$ (recall Proposition \ref{iper-reduced G-degree(d=4)} and Theorem \ref{Th.1}).

Moreover, if $(\Gamma,\gamma)$ represents a singular $4$-manifold $M^4$,  Proposition \ref{iper-reduced G-degree(d=4)} and Proposition \ref{inequalities Euler characteristic} imply:
$$ \sum_{i\in \Delta_{4}} \omega_G(\Gamma_{\hat i}) \  = \  3 \Big(2 \chi(M^4) + p - \sum_{i\in \Delta_{4}}  g_{\hat i}\Big).$$
Note that, if $(\Gamma,\gamma)$ is a crystallization of a closed PL $4$-manifold $M^4$ with $rk(\pi_1(M^4))=m$, relation \eqref{p and q} gives:
$$\sum_{i\in \Delta_{4}} \omega_G(\Gamma_{\hat i}) \  =
\ 3 \Big[5 \Big(\chi(M^4) + 2m -2 \Big) + q \Big].$$
In particular, if $(\Gamma,\gamma)$ is semi-simple (i.e. $q=0$), $ \sum_{i\in \Delta_{4}} \omega_G(\Gamma_{\hat i}) \ = \ 15 \Big(\chi(M^4) + 2m -2\Big)$ follows, as \cite[Proposition 8]{[Basak-Casali 2016]} trivially implies.
}
\end{remark}

\medskip\medskip\medskip

\section*{Acknowledgments}
This work was supported by the {\it ``National Group for Algebraic and Geometric Structures, and their Applications''} (GNSAGA - INDAM) and by University of Modena and Reggio Emilia, projects  {\it ``Colored graphs representing pseudomanifolds: an interaction with random geometry and physics"} and {\it ``Applicazioni della Teoria dei Grafi nelle Scienze, nell'Industria e nella Societ\'a"}.

\smallskip
\par \noindent
The authors would like to thank Gloria Rinaldi for her helpful ideas and suggestions about relationship between cyclic permutation properties and Hamiltonian cycle decompositions of complete graphs.

\noindent They also thank the referees for their very helpful suggestions.

\end{document}